\documentclass[11pt]{amsart}
\usepackage{amsmath}
\usepackage{amssymb}
\usepackage{amscd}

\def\be{\begin{equation}}       \def\ee{\end{equation}}
\def\bd{\begin{displaymath}}    \def\ed{\end{displaymath}}
\def\beq{\begin{eqnarray}}      \def\eeq{\end{eqnarray}}
\def\bseq{\begin{eqnarray*}}    \def\eseq{\end{eqnarray*}}
\def\ba{\begin{array}}          \def\ea{\end{array}}
\def\ben{\begin{enumerate}}     \def\een{\end{enumerate}}

\def\ra{\rightarrow}

\def\lra{\longrightarrow}

\def\cop{\Delta}

\def\cnt{\varepsilon}
\def\ot{\otimes}

\def\GLqtwo{{GL}_{q}(2)}
\def\GLhtwo{{GL}_{h}(2)}

\def\GLhh'two{{GL}_{h,h'}(2)}

\def\Grs{{G}_{r,s}}

\def\Grss'{{G}_{r}^{s,s'}}
\def\Gmkk'{{G}_{m}^{k,k'}}

\def\ident{{\bf 1}}

\def\ca{{\mathcal A}}

\begin{document}

\title{THE COLOURED QUANTUM PLANE}

\author{Deepak Parashar}

\address{Max-Planck-Institute for Mathematics in the Sciences \\
Inselstrasse 22-26 \\ D-04103 Leipzig \\ Germany}

\email{Deepak.Parashar@mis.mpg.de}
\urladdr{http://personal-homepages.mis.mpg.de/parashar/}
\keywords{}

\subjclass{}

\begin{abstract}

We study the quantum plane associated to the coloured quantum group
$GL_{q}^{\lambda,\mu}(2)$ and solve the problem of constructing the
corresponding differential geometric structure. This is achieved within
the $R$-matrix framework generalising the Wess-Zumino formalism and leads 
to the concept of coloured quantum space. Both, the coloured Manin plane
as well as the bicovariant differential calculus exhibit the colour
exchange symmetry. The coloured $h$-plane corresponding to the
coloured Jordanian quantum group $GL_{h}^{\lambda,\mu}(2)$  is also
obtained by contraction of the coloured $q$-plane. \\ 
{[}{\it Jour. Geom. Phys.}, to appear.] \\
{\it MSC (2000):} 81R50; 17B37 \\
{\it Keywords:} Quantum groups; Manin plane; Differential calculus
\end{abstract}
\maketitle

\section{Introduction}

Manin's approach \cite{manin} to quantum groups is based on the fact that
these structures can be considered as comodule algebras for the so-called
{\it quantum} or {\it Manin} plane. In analogy with the classical case, a
quantum group also acts upon a formal vector space called the quantum
vector space, or simply {\it quantum space}. For the simplest example of
the two-dimensional case, this is generated by elements $x$ and $y$
satisfying the nontrivial commutation relation
\begin{equation}
xy=qyx
\label{qplane}
\end{equation}
where $q\neq1$ is the quantum deformation parameter and $\{ x,y \}$
commute with the generating elements of the quantum group $\GLqtwo$. In
the past decade or so, the study of quantum planes and quantum spaces has
been at the forefront of mathematical physics research. One of the most
interesting aspects is the formulation of a consistent differential
calculus on the quantum plane, a powerful procedure given by Wess and
Zumino \cite{wess}. An immediate consequence is the $q$-deformation of the
quantum mechanical phase space, also given in \cite{wess}. Since then,
extensive effort has been directed towards investigating the
algebraic and differential geometric structure of quantum planes and
quantum spaces in the noncommutative setting such as development of
noncommutative gauge theories and quantum field theories. The study of
such geometric structures also lies at the interface of Connes'
noncommutative differential geometry and the theory of quantum groups. The
Wess-Zumino formalism admits interesting multiparametric generalisations
\cite{thesis} and the technique has also been successfully applied to the
$h$-planes \cite{agha,madore,book} associated with the nonstandard
$h$-deformations. 
\\
More recently, `coloured' quantum groups have been studied
\cite{basu,quesne,preeti,dpjmp,zhang} in the context of generalisation of
the basic theory of quantum groups by parametrising the corresponding
generators by some continuously varying colour parameters. Emerged
originally as non-additive solutions of the spectral parameter dependent
quantum Yang-Baxter equation in the study of integrable systems, coloured
quantum groups are infinite-dimensional algebras. This is an interesting
feature since on the one hand these are examples of infinite-dimensional
Hopf algebras, now known as coloured Hopf algebras, while on the other
hand they can be couched in the language familiar to those for the
finite-dimensional ones. These are also related to invariants of knots
and braid group respresentations similar to the usual uncoloured quantum
groups. For further motivation see \cite{quesne,zhang} and references
therein.
\\
Despite the interest generated in the study of the coloured extension of
quantum groups, little is known about the quantum planes or
quantum spaces associated to these structures and the whole theory remains
pretty much in its embryonic stage. In a recent work \cite{dplmp}, the
coloured quantum group  $GL_{q}^{\lambda,\mu}(2)$ was investigated to
establish a duality between the coloured algebra of quantised functions on
the group and that of its universal enveloping algebra, i.e., its dual. In
\cite{dplmp}, the author also gave a coloured generalisation of the
$R$-matrix approach to construct a bicovariant differential calculus on 
$GL_{q}^{\lambda,\mu}(2)$. A natural question that arises at this stage
pertains to the quantum spaces upon which coloured quantum groups act. 
It is the purpose of this paper, therefore, to present a differential
calculus on the coloured quantum planes associated to 
$GL_{q}^{\lambda,\mu}(2)$. This is carried out by generalising the
Wess-Zumino formalism and leads to the concept of a {\it coloured
quantum space}.
\section{Coloured $FRT$ algebra}
It is well-established \cite{basu} that $GL_{q}^{\lambda,\mu}(2)$ provides
a coloured generalisation of $\GLqtwo$ resulting in a coloured
Faddeev-Reshetikhin-Takhtajan $(FRT)$ \cite{frt} algebra. As already
mentioned, the corresponding generators are parametrised by
continuously varying colour parameters $\lambda$ and $\mu$, and as such
the whole algebra becomes infinite-dimensional. Recall that the
$GL_{q}^{\lambda,\mu}(2)$ commutation algebra is given as
\begin{equation}
\begin{array}{ll}
a_{\lambda}b_{\mu}=q^{1+2\lambda}b_{\mu}a_{\lambda} \qquad &
a_{\lambda}c_{\mu}=q^{1-2\lambda}c_{\mu}a_{\lambda}\\
b_{\lambda}c_{\mu}=q^{-2(\lambda+\mu)}c_{\mu}b_{\lambda} \qquad &
b_{\lambda}d_{\mu}=q^{1-2\mu}d_{\mu}b_{\lambda}\\
d_{\lambda}c_{\mu}=q^{-1-2\lambda}c_{\mu}d_{\lambda} \qquad &
\left[a_{\lambda},d_{\mu}\right]=(q-q^{-1})q^{\lambda+\mu}b_{\mu}c_{\lambda}
\end{array}
\end{equation}
\begin{equation}
\begin{array}{ll} 
a_{\lambda}b_{\mu}=q^{\lambda-\mu}a_{\mu}b_{\lambda} \qquad &
a_{\lambda}c_{\mu}=q^{\mu-\lambda}a_{\mu}c_{\lambda}\\
b_{\lambda}c_{\mu}=b_{\mu}c_{\lambda} \qquad &
b_{\lambda}d_{\mu}=q^{\lambda-\mu}b_{\mu}d_{\lambda}\\
c_{\lambda}d_{\mu}=q^{\mu-\lambda}c_{\mu}d_{\lambda} \qquad &
a_{\lambda}d_{\mu}=a_{\mu}d_{\lambda}
\end{array}
\end{equation}
\begin{equation}
\begin{array}{lll}
a_{\lambda}a_{\mu} &=& a_{\mu}a_{\lambda}\\
b_{\lambda}b_{\mu} &=& q^{2(\lambda-\mu)}b_{\mu}b_{\lambda}\\
c_{\lambda}c_{\mu} &=& q^{2(\mu-\lambda)}c_{\mu}c_{\lambda}\\
d_{\lambda}d_{\mu} &=& d_{\mu}d_{\lambda}
\end{array}
\end{equation}
which is obtained from the coloured $RTT$ relations where
\begin{equation}
R(\lambda,\mu) = \begin{pmatrix}
q^{1-\lambda+\mu} & 0 & 0 & 0\\
0 & q^{\lambda+\mu} & 0 & 0\\
0 & q-q^{-1} & q^{-\lambda-\mu} & 0\\
0 & 0 & 0 & q^{1+\lambda-\mu}\
\end{pmatrix}
\end{equation}
is non-additive and satisfies the coloured quantum Yang-Baxter equation,
and
$T_{\lambda}= \left(
\begin{smallmatrix}a_{\lambda}&b_{\lambda}\\c_{\lambda}&d_{\lambda}
\end{smallmatrix}\right)$ and
$T_{\mu}= \left(
\begin{smallmatrix}a_{\mu}&b_{\mu}\\c_{\mu}&d_{\mu}\end{smallmatrix} 
\right)$
are quantum matrices of the generators. The coalgebra structure is given
by
$\cop (T_{\lambda})=T_{\lambda}\dot{\ot} T_{\lambda}$,
$\cnt (T_{\lambda})=\ident$. 
The quantum determinant is
$D_{\lambda}=a_{\lambda}d_{\lambda}-q^{1-2\lambda}c_{\lambda}b_{\lambda}$
and the antipode is
\begin{equation}
S(T_{\lambda})= 
D_{\lambda}^{-1}\begin{pmatrix}d_{\lambda}&-q^{-1+2\lambda}b_{\lambda}\\
-q^{1-2\lambda}c_{\lambda}&a_{\lambda}\end{pmatrix}
\end{equation}
Clearly, the colourless $\GLqtwo$ $FRT$ algebra is retrieved with
vanishing of the colour parameters.
\section{Coloured quantum plane}
In analogy with the standard uncoloured quantum group $\GLqtwo$,
its coloured version $GL_{q}^{\lambda,\mu}(2)$ can be considered as an
endomorphism of a pair of quadratic algebras $\ca_{q}^{2|0}$ and
$\ca_{q}^{0|2}$ (its dual) in the sense of Manin \cite{manin}. The algebra
$\ca_{q}^{2|0}$ is generated by elements $\{ x_{\lambda}, x_{\mu},
y_{\lambda}, y_{\mu} \}$, the `coloured' coordinates, such that
\begin{equation}
\begin{array}{lll}
x_{\lambda}x_{\mu} &=& q^{\lambda-\mu}x_{\mu}x_{\lambda}\\
x_{\lambda}y_{\mu} &=& q^{1-\lambda-\mu}y_{\mu}x_{\lambda}\\
x_{\lambda}y_{\mu} &=& x_{\mu}y_{\lambda}\\
y_{\lambda}y_{\mu} &=& q^{\mu-\lambda}y_{\mu}y_{\lambda}
\end{array}
\label{boson}
\end{equation}
This is the so-called {\it coloured quantum plane} suggested in
\cite{basu} and defines a left coaction
\begin{equation}
\cop_{L}: \ca_{q}^{2|0} \lra GL_{q}^{\lambda,\mu}(2) \ot \ca_{q}^{2|0}
\end{equation}
i.e.
\begin{equation}
\begin{array}{lll}
\cop_{L}\left( 
\begin{smallmatrix}x_{\lambda}\\y_{\lambda}\end{smallmatrix}\right)
&\lra& \left( \begin{smallmatrix} 
a_{\lambda}&b_{\lambda}\\c_{\lambda}&d_{\lambda}\end{smallmatrix}\right) 
\ot \left(
\begin{smallmatrix}x_{\lambda}\\y_{\lambda}\end{smallmatrix}\right) \\
\cop_{L}\left( \begin{smallmatrix}x_{\mu}\\y_{\mu}\end{smallmatrix}\right)
&\lra& \left( \begin{smallmatrix}
a_{\mu}&b_{\mu}\\c_{\mu}&d_{\mu}\end{smallmatrix}\right) \ot
\left( \begin{smallmatrix}x_{\mu}\\y_{\mu}\end{smallmatrix}\right)
\end{array}
\end{equation}
where the coloured quantum matrices $T_{\lambda}$ and $T_{\mu}$
commute with the coloured coordinates $\{ x_{\lambda}, x_{\mu},
y_{\lambda}, y_{\mu} \}$. In general, this can be written as
\begin{equation}
x^{i}_{\lambda}x^{j}_{\mu}=B^{ij}_{kl}x^{k}_{\mu}x^{l}_{\lambda}
\label{interior}
\end{equation}
where
\begin{equation}
B^{ij}_{kl}=\left \{
\begin{array}{l}B(\lambda,\mu)=\frac{1}{q}\hat{R}(\lambda,\mu);\qquad i<j\\
B(\mu,\lambda)=\frac{1}{q}\hat{R}(\mu,\lambda); \qquad i=j
\end{array}\right.
\end{equation}
and
\begin{equation}
\hat{R}(\lambda,\mu) = \begin{pmatrix}
q^{1-\lambda+\mu} & 0 & 0 & 0\\
0 & q-q^{-1} & q^{-\lambda-\mu} & 0\\
0 & q^{\lambda+\mu} & 0 & 0\\
0 & 0 & 0 & q^{1+\lambda-\mu}\
\end{pmatrix}
\end{equation}
is the coloured braided $R$-matrix which solves the coloured
braided Yang-Baxter equation. The right coaction is given by
$\ca_{q}^{0|2}$, the algebra dual to $\ca_{q}^{2|0}$
\begin{equation}
\cop_{R}: \ca_{q}^{0|2} \lra \ca_{q}^{0|2} \ot GL_{q}^{\lambda,\mu}(2)
\end{equation}
i.e.
\begin{equation}
\begin{array}{lll}
\cop_{R}\left(
\begin{smallmatrix}\xi_{\lambda}&\eta_{\lambda}\end{smallmatrix}\right)
&\lra& \left(
\begin{smallmatrix}\xi_{\lambda}&\eta_{\lambda}\end{smallmatrix}\right) 
\ot \left( \begin{smallmatrix}
a_{\lambda}&b_{\lambda}\\c_{\lambda}&d_{\lambda}\end{smallmatrix}\right)\\
\cop_{R}\left(
\begin{smallmatrix}\xi_{\mu}&\eta_{\mu}\end{smallmatrix}\right)
&\lra& \left(
\begin{smallmatrix}\xi_{\mu}&\eta_{\mu}\end{smallmatrix}\right)
\ot \left( \begin{smallmatrix}
a_{\mu}&b_{\mu}\\c_{\mu}&d_{\mu}\end{smallmatrix}\right)
\end{array}   
\end{equation}
The algebra $\ca_{q}^{0|2}$ is generated by the coloured Grassmann
variables $\{ \xi_{\lambda},\xi_{\mu},\eta_{\lambda},\eta_{\mu} \}$
satisfying
\begin{equation}
\begin{array}{lll}
\xi_{\lambda}\xi_{\mu} &=& \eta_{\lambda}\eta_{\mu} = 0\\
\xi_{\lambda}\eta_{\mu} &=& -q^{-(1+\lambda+\mu)}\eta_{\mu}\xi_{\lambda}\\
\xi_{\lambda}\eta_{\mu} &=& \xi_{\mu}\eta_{\lambda}
\end{array}
\label{fermion}
\end{equation}
and defines the {\it coloured quantum hyperplane}. This can also be
expressed in a more general form
\begin{equation}
\xi^{i}_{\lambda}\xi^{j}_{\mu}=-C^{ij}_{kl}\xi^{k}_{\mu}\xi^{l}_{\lambda}
\label{exterior}
\end{equation}
where
\begin{equation}
C^{ij}_{kl}=\left \{
\begin{array}{l}C(\lambda,\mu)=q\hat{R}(\lambda,\mu);\qquad i<j\\
C(\mu,\lambda)=q\hat{R}(\mu,\lambda);\qquad i=j
\end{array}\right.
\label{cmatrix}
\end{equation}
Both sets of relations (\ref{boson}) as well as (\ref{fermion}) of the
coloured quantum plane satisfy the $\lambda \leftrightarrow \mu$ exchange
symmetry and are invariant under the coaction of 
$GL_{q}^{\lambda,\mu}(2)$. In other words, the coloured quantum group
$GL_{q}^{\lambda,\mu}(2)$ acts on a quantum vector space generated by the
coloured variables $\{ x_{\lambda}, x_{\mu},y_{\lambda}, y_{\mu} \}$ and
this defines a `coloured quantum space'. In the limit of the vanishing
colour parameters, we recover exactly the uncoloured two-dimensional
quantum plane (\ref{qplane}) corresponding to $\GLqtwo$.
\section{Differential calculus}
Before we proceed to the construction of a differential
calculus on the coloured quantum plane, it is important to stress that
equations (\ref{interior}) and (\ref{exterior}) are coloured 
generalisations of the usual
relations for the Manin plane and its exterior algebra, and so are the 
matrices $B^{ij}_{kl}$ and $C^{ij}_{kl}$. Therefore, the Wess-Zumino
calculus \cite{wess} has to be modified in order to incorporate colour
dependence. The resulting calculus is much more complicated, has colour
copies of the various mathematical quantities and for convenience we
follow the same notation as \cite{wess}. For the derivatives, we define
\begin{eqnarray}
\partial_{x_{\lambda}}=\frac{\partial}{\partial_{x_{\lambda}}},\qquad
\partial_{x_{\mu}}=\frac{\partial}{\partial_{x_{\mu}}}\\
\partial_{y_{\lambda}}=\frac{\partial}{\partial_{y_{\lambda}}},\qquad
\partial_{y_{\mu}}=\frac{\partial}{\partial_{y_{\mu}}}
\end{eqnarray}
and for the differentials
\begin{eqnarray}
\xi_{\lambda}=\mathsf{d}x_{\lambda}, \qquad \xi_{\mu}=\mathsf{d}x_{\mu}\\
\eta_{\lambda}=\mathsf{d}y_{\lambda}, \qquad \eta_{\mu}=\mathsf{d}y_{\mu}
\end{eqnarray}
The differential calculus on the coloured quantum plane $\ca_{q}^{2|0}$ is
given by the commutation relations between the variables and the
derivatives
\begin{equation}
\partial_{j_{\lambda}}x^{i}_{\mu}=\delta^{i}_{j} +
C^{ik}_{jl}x^{l}_{\mu}\partial_{k_{\lambda}}
\end{equation}
with matrix $C$ defined in (\ref{cmatrix}). For $i,j=1..2$, this yields
for $\ca_{q}^{2|0}$
\begin{equation}
\begin{array}{lll}
\partial_{x_{\lambda}}y_{\mu} &=&
q^{1+\lambda+\mu}y_{\mu}\partial_{x_{\lambda}}\\
\partial_{y_{\lambda}}x_{\mu} &=& 
q^{1-\lambda-\mu}x_{\mu}\partial_{y_{\lambda}}\\
\partial_{y_{\lambda}}y_{\mu} &=&
1+q^{2-\lambda+\mu}y_{\mu}\partial_{y_{\lambda}}\\
\partial_{x_{\lambda}}x_{\mu} &=&
1+q^{2+\lambda-\mu}x_{\mu}\partial_{x_{\lambda}} +
(q^{2}-1)y_{\mu}\partial_{y_{\lambda}}
\end{array}
\end{equation}
In addition, the derivatives also satisfy relations among themselves
\begin{equation}
\partial_{i_{\lambda}}\partial_{j_{\mu}} =
F^{lk}_{ji}\partial_{k_{\mu}}\partial_{l_{\lambda}};\qquad i\leq j
\end{equation}
where
\begin{equation}
F^{lk}_{ji}=F(\lambda,\mu)=\frac{1}{q}\hat{R}(\lambda,\mu)
\end{equation}
This gives
\begin{equation}
\partial_{x_{\lambda}}\partial_{y_{\mu}} =
q^{-(1+\lambda+\mu)}\partial_{y_{\mu}}\partial_{x_{\lambda}}
\end{equation}
The differential calculus on the coloured quantum hyperplane
$\ca_{q}^{0|2}$ is given by the commutation relations between the
variables and the differentials
\begin{equation}
x^{i}_{\lambda}\xi^{j}_{\mu}=C^{ij}_{kl}\xi^{k}_{\mu}x^{l}_{\lambda}
\end{equation}
Again, the matrix $C$ has the same form as (\ref{cmatrix}). Explicitly,
for $\ca_{q}^{0|2}$ we obtain
\begin{equation}
\begin{array}{lll}
x_{\lambda}\xi_{\mu} &=& q^{2+\lambda-\mu}\xi_{\mu}x_{\lambda}\\
y_{\lambda}\xi_{\mu} &=& q^{1+\lambda+\mu}\xi_{\mu}y_{\lambda}\\
y_{\lambda}\eta_{\mu} &=& q^{2-\lambda+\mu}\eta_{\mu}y_{\lambda}\\
x_{\lambda}\eta_{\mu} &=& q^{1-\lambda-\mu}\eta_{\mu}x_{\lambda} +
(q^{2}-1)\xi_{\mu}y_{\lambda}
\end{array}
\end{equation}
Furthermore, the commutation relations between the derivatives and the
differentials are
\begin{equation}
\partial_{j_{\lambda}}\xi^{i}_{\mu}=
D^{ik}_{jl}\xi^{l}_{\mu}\partial_{k_{\lambda}}
\end{equation}
where
\begin{equation}
D^{ik}_{jl}=\left \{
\begin{array}{l}D(\lambda,\mu)=C^{-1}(\lambda,\mu);\qquad i<j\\ 
D(\mu,\lambda)=C^{-1}(\mu,\lambda);\qquad i=j
\end{array}\right.
\end{equation}  
So, we have for $\ca_{q}^{0|2}$
\begin{equation}
\begin{array}{lll}
\partial_{x_{\lambda}}\xi_{\mu} &=&
\frac{1}{q^{2+\lambda-\mu}}\xi_{\mu}\partial_{x_{\lambda}}\\
\partial_{x_{\lambda}}\eta_{\mu} &=& 
\frac{1}{q^{1-\lambda-\mu}}\eta_{\mu}\partial_{x_{\lambda}}\\
\partial_{y_{\lambda}}\xi_{\mu} &=&
\frac{1}{q^{1+\lambda+\mu}}\xi_{\mu}\partial_{y_{\lambda}}\\
\partial_{y_{\lambda}}\eta_{\mu} &=& 
\frac{1}{q^{2-\lambda+\mu}}\eta_{\mu}\partial_{y_{\lambda}} +
(\frac{1}{q^{2}}-1)\xi_{\mu}\partial_{x_{\lambda}}
\end{array}
\end{equation}
Finally, the exterior differential is defined as
\begin{equation}
\begin{array}{c}
\mathbf{d} = \xi^{i}_{\lambda}\partial_{i_{\lambda}} = 
\xi^{i}_{\mu}\partial_{i_{\mu}}\\
\mathbf{d} = \xi_{\lambda}\partial_{x_{\lambda}} +
\eta_{\lambda}\partial_{y_{\lambda}} = \xi_{\mu}\partial_{x_{\mu}}
+ \eta_{\mu}\partial_{y_{\mu}}
\end{array}
\end{equation}
It can be checked that the exterior differential is nilpotent
($\mathbf{d}^{2}=0$) and satisfies the Leibniz rule
$\mathbf{d}(fg)=(\mathbf{d}f)g+f(\mathbf{d}g)$. The full differential
calculus on the coloured quantum plane and the hyperplane respects the
$\lambda \leftrightarrow \mu$ exchange symmetry, satisfies all consistency
conditions and is invariant under the coaction of the coloured quantum
group $GL_{q}^{\lambda,\mu}(2)$. Note that for vanishing colour
parameters, the calculus reduces exactly to the calculus on the standard
uncoloured two-dimensional quantum plane corresponding to $\GLqtwo$. It is
also important to note that while we have treated the coloured version of
the two-dimensional quantum plane, the result is actually 
infinite-dimensional since $\lambda$ and $\mu$ are continuously
varying colour parameters. As such, not only are the coloured quantum
planes infinite-dimensional but the differential calculus obtained
here is also infinite-dimensional.
\section{Coloured $h$-plane}
We now look at the nonstandard (or Jordanian) counterpart of the
$q$-plane, known as the $h$-plane \cite{agha,madore}. Motivated by the
observation that coloured Jordanian quantum groups can be obtained from
their $q$-deformed counterparts \cite{agha, dpjmp, agha2} using the
well-known contraction procedure, we perform the contraction
on the coloured $q$-plane discussed above. We consider the transformation
on the coloured coordinates
\begin{equation}
\left(
\begin{smallmatrix}x_{\lambda}\\y_{\lambda}\end{smallmatrix}\right) \lra 
g\left( 
\begin{smallmatrix}x_{\lambda}\\y_{\lambda}\end{smallmatrix}\right) \qquad
\text{and} \qquad
\left( 
\begin{smallmatrix}x_{\mu}\\y_{\mu}\end{smallmatrix}\right) \lra  
g\left(
\begin{smallmatrix}x_{\mu}\\y_{\mu}\end{smallmatrix}\right)
\end{equation}
where $g$ is the two-dimensional transformation matrix
\begin{equation}
g=\left( \begin{smallmatrix}1&\alpha\\0&1\end{smallmatrix} \right)
\end{equation}
and $\alpha=\frac{h}{q-1}$ depends on both `$q$' and the `$h$' deformation
parameters. Applying this transformation to the coloured $q$-plane
(\ref{boson}), we obtain in the singular limit $q\ra 1$
\begin{equation}
[x_{\lambda}, y_{\mu}] = h(1-2\mu)y_{\lambda}y_{\mu}
\end{equation}
i.e. the coloured $h$-plane corresponding to the coloured Jordanian
quantum group $GL_{h}^{\lambda,\mu}(2)$ \cite{preeti,dpjmp}. Again, for
vanishing colour parameters, this reduces to the well-known $h$-plane
$[x,y]=hy^{2}$ for the Jordanian quantum group $\GLhtwo$. While performing
the above transformation, we also obtain a `hybrid' standard-nonstandard
quantum plane
\begin{equation}
x_{\lambda}y_{\mu}-q^{1-\lambda-\mu}y_{\mu}x_{\lambda} =
h[1-2\mu]_{q}y_{\lambda}y_{\mu} =
hq^{\mu-\lambda}[1-2\mu]_{q}y_{\mu}y_{\lambda} 
\end{equation}
which we call the `coloured $(q,h)$-plane', where
$[x]_{q}=\frac{1-q^{x}}{1-q}$ is the basic number from $q$-analysis. 
\section{Conclusions}
We have presented the coloured quantum plane and the hyperplane
covariant under the action of the coloured quantum group
$GL_{q}^{\lambda,\mu}(2)$. Generalising the Wess-Zumino formalism, we have
also constructed a consistent differential calculus on the coloured
quantum plane within the framework of noncommutative geometry.  In
addition, we have performed a contraction on the coloured $q$-plane to
obtain the coloured $h$-plane corresponding to the Jordanian quantum group 
$GL_{h}^{\lambda,\mu}(2)$. The results obtained are valid also for
higher-dimensional and multiparametric coloured quantum planes.
\\
Before closing, it is pertinent to comment on a few applications of this
work. Firstly, the calculus developed here equips us with essential
ingredients necessary to develop further the notion of coloured quantum
spaces vis-a-vis coloured version of the $q$-deformed quantum mechanical
phase space. Secondly, it would be useful to construct a differential
calculus on the coloured $h$-plane which could possibly be obtained by
contraction of the calculus on the coloured $q$-plane presented here, and
to investigate the intermediate coloured $(q,h)$-planes. One could then
build up a whole class of differential calculi on various different
coloured quantum planes including the supersymmetric
(graded) versions. The coloured Manin plane also forms interesting example
of quantum homogeneous spaces \cite{tarlini}, the geometric nature of
which is currently being studied from several aspects.
\\
Besides, it will be useful to investigate differential
geometric structure on the underlying coloured Hopf algebras, which
arose in the knot theoretic context. Also important to address are issues
such as ring theoretic and quantum matrix theory aspects of coloured
quantum groups. It would be of interest to explore further relation
with the already well-known $q$-Heisenberg algebras, as well as the
possibility of coloured generalisation of some new  quantum
groups \cite{taft} given recently.


\begin{thebibliography}{99}
\bibitem{manin}
Yu. I. Manin, Quantum groups and noncommutative geometry, {\it Preprint
Montreal University} CRM-1561 (1988).
\bibitem{wess}
J. Wess and B. Zumino, Covariant differential calculus on the quantum
hyperplane, {\it Nucl. Phys. (Proc. Suppl.)} {\bf 18B} (1990) 302 - 312.
\bibitem{thesis}
P. Parashar, {\it Multiparametric Deformations of Quantum Vector
Spaces}, Ph.D. Thesis, University of Delhi, 1992.
\bibitem{agha}
A. Aghamohammadi, The two-parametric extension of $h$-deformation of
$GL(2)$, and the differential calculus on its quantum plane, {\it Modern
Phys. Lett. A} {\bf 8} (1993) 2607 - 2613. 
\bibitem{madore}
S. Cho, J. Madore and K. S. Park, Noncommutative geometry of the
$h$-deformed quantum plane, {\it J. Phys. A} {\bf 31} (1998) 2639 - 
2654.
\bibitem{book}
J. Madore, {\it An Introduction to Noncommutative Differential
Geometry and its Physical Applications}, Cambridge University Press, 1995.
\bibitem{basu}
B. Basu-Mallick, Coloured extensions of $\GLqtwo$ quantum group and
related noncommutative planes, {\it Internat. J. Modern Phys. A} {\bf 10},
(1995) 2851 - 2864.
\bibitem{quesne}
C. Quesne, Coloured quantum universal enveloping algebras, {\it
J. Math. Phys.} {\bf 38} (1997) 6018 - 6039; Duals of coloured
quantum universal enveloping algebras and coloured universal $T$-matrices,
{\it J. Math. Phys.} {\bf 39} (1998) 1199 - 1222.
\bibitem{preeti}
P. Parashar, Jordanian $U_{h,s}gl(2)$ and its coloured realisation, {\it
Lett. Math. Phys.} {\bf 45} (1998) 105 - 112.
\bibitem{dpjmp}
D. Parashar and R. J. McDermott, Contraction of the $\Grs$ quantum
group to its nonstandard analogue and corresponding coloured quantum
groups, {\it J. Math. Phys.} {\bf 41} (2000) 2403 - 2416.
\bibitem{zhang}
P.-G. Luan, H. C. Lee and R. B. Zhang, Coloured solutions of the
Yang-Baxter equation from representations of $U_{q}gl(2)$, {\it
J. Math. Phys.} {\bf 41} (2000) 6529 - 6543.
\bibitem{dplmp}
D. Parashar, Duality for coloured quantum groups, {\it Lett. Math. Phys.}
{\bf 53} (2000) 29 - 40.
\bibitem{frt}
L. D. Faddeev, N. Y. Reshetikhin and L. A. Takhtajan, Quantisation of
Lie groups and Lie algebras, {\it Leningrad Math. J.} {\bf 1} (1990) 193 -
225.
\bibitem{agha2}
A. Aghamohammadi, M. Khorrami and A. Shariati, Jordanian deformation of
SL(2) as a contraction of its Drinfeld-Jimbo deformation, {\it J. Phys. A}
{\bf 28} (1995) L225 - L331.
\bibitem{tarlini}
F. Bonechi, N. Ciccoli, R. Giachetti, E. Sorace and M. Tarlini, The
coisotropic subgroup structure of $SL_{q}(2,\mathbb{R})$, {\it
J. Geom. Phys.} {\bf 37} (2001) 190 - 200, and references therein.
\bibitem{taft}
S. Rodriguez-Romo and E. Taft, private communication (2001).
\end{thebibliography}
\end{document}